# On the intermediate asymptotic efficiency of goodness-of-fit tests in multinomial distributions

Sherzod M. Mirakhmedov

V.I. Romanovskiy Institute of Mathematics. Academy of Sciences of Uzbekistan
University str., 9, Tashkent -100174, Uzbekistan
e-mail: shmirakhmedov@yahoo.com

**Abstract**. We consider goodness-of-fit tests for uniformity of a multinomial distribution by means of tests based on a class of symmetric statistics, defined as the sum of some function of cell-frequencies. We are dealing with an asymptotic regime, where the number of cells grows with the sample size. Most attention is focused on the class of power divergence statistics. The aim of this article is to study the intermediate asymptotic relative efficiency of two tests, where the powers of the tests are asymptotically non-degenerate and the sequences of alternatives converge to the hypothesis, but not too fast. The intermediate asymptotic relative efficiency of the chi-square test wrt an arbitrary symmetric test is considered in details.

**Keywords.** Asymptotic efficiency, chi-square statistic, log-likelihood ratio statistic, goodness-of-fit tests, multinomial distribution, power divergence statistics.

**MSC (2000):** 62G10, 62G20

## 1, Introduction

Let $(\eta_1,...,\eta_N)$ be a random vector of frequencies of a multinomial model on $n \geq 1$ observations classified into $N > 1$ cells with the cell-probabilities $P = (p_1,...,p_N)$, $p_1 + ... + p_N = 1$, all $p_j > 0$. We assume that $N = N(n) \to \infty$ as $n \to \infty$, so that $n/N \to \lambda \in [0, \infty]$ and $n^2/N \to \infty$. We are interested in the problem of testing of $H_0$: $P = (N^{-1},...,N^{-1})$ versus sequences of alternatives $H_{1n}$: $P = (p_1,...,p_N) \neq (N^{-1},...,N^{-1})$, which *approach* $H_0$ as $n \to \infty$ so that

$$\varepsilon(n) = \frac{1}{N}\sum_{m=1}^{N}\varepsilon_{m,n}^2 \to 0, \qquad (1.1)$$

where $\varepsilon_{m,n} = Np_m - 1$, i.e., $p_m = N^{-1}(1+\varepsilon_{m,n})$, $m = 1,...,N$ and $\max_m |\varepsilon_{m,n}| \to 0$ as $n \to \infty$. We consider a class of symmetric tests based on statistics of the form

$$S_n^h = \sum_{l=1}^{N} h(\eta_l), \qquad (1.2)$$

where $h$ is a nonlinear real-valued function, defined on the non-negative axis. A test based on the statistic $S_n^h$, whose large values reject $H_0$, is called $h$-test for short.

The *power divergence statistics* (PDS) $CR_n(d)$ of Cressie and Read (1984) is an important sub-class of statistics (1.2), where $h(x) = \psi_d(x)$, $d > -1$,

$$\psi_d(x) = 2(d(d+1))^{-1} x[(x/\lambda_n)^d - 1], \ d \neq 0, \text{ else } \psi_0(x) = 2x \log(x/\lambda_n), \quad (1.3)$$

$\lambda_n = n/N$ is the average of observations arisen in the cells. (One could consider $d \to -1$ and $d < -1$, but this unnecessary in the context of this paper). The most commonly used special versions of PDS are the chi-square statistic $\chi_n^2 = CR_n(1)$, the log-likelihood ratio statistic $\Lambda_n = CR_n(0)$, and the Freeman-Tukey statistic $T_n^2 = CR_n(-1/2)$, viz.,

$$\chi_n^2 = \lambda_n^{-1} \sum_{m=1}^{N} (\eta_m - \lambda_n)^2, \ \Lambda_n = 2 \sum_{m=1}^{N} \eta_m \log(\eta_m / \lambda_n), \ T_n^2 = 4 \sum_{m=1}^{N} (\sqrt{\eta_m} - \sqrt{\lambda_n})^2, \quad (1.4)$$

Count statistics (CS) form another sub-class of (1.2). Two popular variants of CS are $\mu_r$ - the number of cells containing exactly $r$ observations and $C_n$ - the number of collisions (i.e., the number of those observations that fall in a cell that already has an observation in it), viz.,

$$\mu_r = \sum_{m=1}^{N} \text{I}\{\eta_m = r\}, \ r \geq 0, \text{ and } C_n = \sum_{m=1}^{N} (\eta_m - 1) \text{I}\{\eta_m > 1\}, \quad (1.5)$$

where $\text{I}\{\cdot\}$ denotes the indicator function.

Our goal is to study the intermediate asymptotic relative efficiencies (IARE) of two $h$-tests, provided that their powers are asymptotically non-degenerate and the alternatives (1.1) are such that

$$\nabla_n \to \infty \text{ and } \nabla_n = o(\sqrt{N}), \quad (1.6)$$

where $\nabla_n = n\varepsilon(n)/\sqrt{N}$.

These alternatives are somewhere between the alternatives which do not approach the hypothesis (or the rate of approach is too slow so that $\varepsilon(n) \geq \lambda_n^{-1}$) and the alternatives which approaches the hypothesis at the rate so that $\nabla_n = O(1)$, i.e., a rate that ensures that the power of $h$-tests of size $\alpha > 0$ has a limit in $(\alpha, 1)$. In other words, we are interested in comparison of two $h$-tests in the situation, which are intermediate between Bahadur's and Pitman's concepts of asymptotic efficiency (AE), see Nikitin(1995), and hence the name. For the detailed comments on the alternatives (1.6) and brief survey we refer to Mirakhmedov (2021). In that work the comparison of $h$-tests was considered in terms of the so called "$\alpha$ – intermediate AE", according to which the performance of an $h$-test of size $\alpha$ measured by the asymptotic value of

$e_N^\alpha(S_n^h) = -\log P_0\{S_n^h > NA_{1n}(h)\}$, provided that its power is asymptotically non-degenerate. Here and throughout the article $A_{in}(h)$ is an asymptotic value of $N^{-1} E_i S_n^h$ (see below section 2), $P_i$ and $E_i S_n^h$, as usually, denote the probability and expectation, counted under $H_0$ and $H_{1n}$, respectively. Among others it is shown that the intermediate asymptotic properties of $h$-tests depend on the parameter $\lambda_n = n/N$, whose asymptotic behavior (as $n \to \infty$) provides the classification of the multinomial model: the *very sparse model* when $\lambda_n \to 0$, $n\lambda_n \to \infty$, the *sparse model* when $\lambda_n \to \lambda \in (0,\infty)$, and the *dense model* when $\lambda_n \to \infty$.

In the present paper we study the IARE of two *h*-tests, defined as a limit of the ratio of sample sizes which guarantee the same precision for both tests: the same significance level tending to 0 slower than exponentially and the same asymptotically non-degenerate power. We remind that the concept of IARE was introduced by Kallenberg (1983) and then elaborated for testing problems by Inglot (1999). Further development and application of this IARE in testing problems has been presented in several works, see Cmiel et al (2020) and references therein. We refer to the paper by Inglot (2020), where so called path wise intermediate efficiency elaborated in Inglot et al (2019) is demonstrated by applying to study the IARE of Neyman–Pearson test wrt Kolmogorov-Smirnov test in the classical problem of testing for uniformity versus a class of heavy-tailed alternatives. We use the idea of Inglot (2020) for discrete models, namely in the problem of testing for the uniformity of a multinomial distribution.

The rest of the paper is organized as follows. Section 2 provides a general result on IARE of two *h*-tests and its application for IARE of the chi-square test wrt *h*-tests. Proofs of all statements and Conclusions are given in Section 3.

## 2. Main Results

In what follows $\xi \sim Poi(\lambda)$ stands for "a random variable (r.v.) $\xi$ has Poisson distribution with parameter $\lambda > 0$"; $\Phi(u)$ denotes a standard normal distribution function; $c_j$ is a positive constant, may not the same in each its occurrence, $\xi$, $\xi_1,...,\xi_N$ are independent r.v.s such that $\xi \sim Poi(\lambda_n)$, $\lambda_n = n/N$, and $\xi_m \sim Poi(np_m)$. *Without mention, all asymptotic statements are considered as $n \to \infty$, hence $N \to \infty$ as well.* Now on $\Im_{alt}$ stands for the family of alternatives (1.1) such that

$$\frac{n\varepsilon(n)}{\sqrt{N}} \to \infty, \quad \frac{n}{N}\max_m \varepsilon_{mn}^2 \to 0. \qquad (2.1)$$

This is a slightly narrower family of alternatives than (1.6), but we have to deal with it because of Proposition 2.1 below. Set

$$g(\xi) = h(\xi) - Eh(\xi) - r_n(\xi - \lambda_n), \quad r_n = \lambda_n^{-1} \operatorname{cov}(h(\xi), \xi)$$

$$\sigma^2(h) = Var\, g(\xi) = Varh(\xi)\left(1 - corr^2(h(\xi), \xi)\right),$$

$$\rho(h, \lambda_n) = corr\left(h(\xi) - r_n\xi, \xi^2 - (2\lambda_n + 1)\xi\right).$$

**Proposition 2.1.** Let sequences of alternatives $H_{1n}$ such that $n\max_m \varepsilon_{mn}^2 / N \to 0$

(i) If $E|h(\xi)\xi^3| < \infty$ then

$$\kappa_n(h) := \sqrt{N}(A_{1,n}(h) - A_{0,n}(h))/\sigma(h) = \sqrt{n\lambda_n/2}\,\varepsilon(n)\rho(h,\lambda_n)(1+o(1)). \tag{2.2}$$

(ii) If $Eh^2(\xi)\xi^2 < \infty$ then

$$\sigma_{1,n}^2(h) = \sigma^2(h)(1+o(1)), \tag{2.3}$$

where

$$A_{i,n}(h) = N^{-1}\sum_{m=1}^{N} E_i h(\xi_m),\ \sigma_{i,n}^2(h) = N^{-1}\sum_{m=1}^{N} Var_i g(\xi_m). \qquad \square$$

Note that under very general set-up (see, for instance, Proposition 2.1 of Mirakhmedov et al (2014))

$$E_i S_n^h = NA_{i,n}(h)(1+o(1))\ \text{and}\ Var_i S_n^h = N\sigma_{i,n}^2(h)(1+o(1)).$$

The following proposition is the consequence of Theorem 1 by Mirakhmedov (1992).

**Proposition 2.2.** Assume

$$E|g(\xi)|^3 / \sigma^3(h)\sqrt{N} \to 0. \tag{2.4}$$

Then

$$P_i\left\{S_n^h < u\sigma_{i,n}(h)\sqrt{N} + NA_{i,n}(h)\right\} = \Phi(u) + o(1),\ i=0,1,$$

where $\sigma_{0,n}^2(h) = \sigma^2(h)$. $\qquad \square$

We emphasize that the condition (2.4) is fulfilled for the sparse models if $E|h(\xi)|^3 < \infty$, and for the very sparse models if $\Delta^2 h(0) \neq 0$, where $\Delta h(x) = h(x+1) - h(x)$. For instance, arbitrary PDS and CS $C_n, \mu_r, r=0,1,2$ satisfy this condition. But for the dense models condition (2.4) may impose an additional condition to $\lambda_n$. For instance, (2.4) is fulfilled for PDS without additional restriction, while, for example, for CS $\mu_r$, $r \geq 0$ and $C_n$ it imposes condition $\lambda_n - \ln N - r\ln\ln N \to -\infty$ and $\lambda_n - \ln N \to -\infty$, respectively.

Notice that $|\rho(h,\lambda_n)| \leq 1$ and $|\rho(h,\lambda_n)| = 1$ iff $h(x) = x^2$, i.e., for the chi-square statistic. But for the CS (1.5) $\rho(h,\lambda_n) = o(1)$ if $\lambda_n \to \infty$. Meaning of the functional $\rho(h,\lambda_n)$ is clarified by the following

**Lemma 2.1.** If $Eh^2(\xi)\xi < \infty$ then $\rho(h,\lambda_n) = corr_0(S_n^h, \chi_n^2) + o(1)$. $\qquad \square$

It follows from Proposition 2.2 that the power of $h$-test of size $\alpha > 0$ asymptotically equal to $\Phi\left(2^{-1/2}\nabla_n |\rho(h,\lambda_n)| - \omega_\alpha\right)$, where $\omega_\alpha = \Phi^{-1}(1-\alpha)$, $\nabla_n = \sqrt{n\lambda_n}\varepsilon(n)$ (see (1.6)). Hence for the alternatives (1.1) such that $\nabla_n = O(1)$ (i.e., for Pitman alternatives) the chi-square test is asymptotically most powerful (AMP) within the class of *h*-tests. The chi-square test is the only AMP for sparse models. But for very sparse and dense models the chi-squared test is no longer unique AMP, since by virtue of Lemma 2.2 and Lemma 2.3 below, there are test statistics for which $\rho(h, \lambda_n) \to 1$ if $\lambda_n \to 0$ or $\lambda_n \to \infty$.

We remark also that for the Pitman alternatives, i.e., $\nabla_n = O(1)$, using Proposition 2.1 (ii) and Proposition 2.2 one can observe, see Mirakhmedov and Bozorov (2021), that Pitman asymptotic relative efficiency of *h*- test wrt $\psi$ -test is

$$PE(S_n^h, S_n^\psi) = \lim_{n \to \infty} \left( \rho^2(S_n^h, \lambda_n) / \rho^2(S_n^\psi, \lambda_n) \right).$$

**Lemma 2.2**. If $\lambda_n \to 0$ and $\Delta^2 h(0) \neq 0$ then

$$\rho(h, \lambda_n) = 1 - \frac{\lambda_n}{6} \left( \frac{\Delta^3 h(0)}{\Delta^2 h(0)} \right)^2 + O(\lambda_n^2).$$

In particularly

$$\rho(\psi_d, \lambda_n) = 1 - \frac{3(3^d - 2^{d+1} + 1)^2}{8(2^d - 1)^2} \lambda_n + O(\lambda_n^2), d \neq 0,$$

$$\rho(\psi_0, \lambda_n) = 1 - \frac{3}{8} \left( \frac{\ln 3/4}{\ln 2} \right)^2 \lambda_n + O(\lambda_n^2),$$

$$\rho(h, \lambda_n) = 1 - c(r)\lambda_n + O(\lambda_n^2), \text{ if } h(x) = I(x = r), r = 0, 1, 2,$$

where $c(0) = 1/6$, $c(1) = 3/8$ and $c(2) = 3/2$. □

**Lemma 2.3**. If $\lambda_n \to \infty$, then for every PDS with parameter $d > -1$ it holds

$$\rho(\psi_d, \lambda_n) = 1 - \frac{(d-1)^2}{6\lambda_n} + O(\lambda_n^{-2}).$$ □

The functional $\rho(h, \lambda_n)$ plays a key role also in determining the intermediate properties of the *h*-tests. From the aforesaid and Mirakhmedov (2021) it follows that an *h*-test for which $\rho(h, \lambda_n)$ close to 1 is preferable. Throughout the rest of the paper the symbol **S** denotes the class of *h*-tests, for which $\lim |\rho(h, \lambda_n)| > 0$ and conditions of Proposition 2.1 are fulfilled.

**2.1. IARE of two h-tests.**

Set $\tilde{S}_n^h = \left( S_n^h - NA_{0,n}(h) \right) / \sigma_{0,n}(h)\sqrt{N}$, $\tilde{S}_n^\psi = \left( S_n^\psi - NA_{0,n}(\psi) \right) / \sigma_{0,n}(\psi)\sqrt{N}$, $\kappa_n(h)$ and $\kappa_n(\psi)$ are defined as in (2.2).

For each *n* and any fixed $z \in \mathbb{R}$ define the significance level of the *h* -test corresponding to the critical value $z + \kappa_n(h)$:

$$\alpha_n = \alpha_n(z) = P_0 \left\{ \tilde{S}_n^h > z + \kappa_n(h) \right\}. \tag{2.5}$$

Due to Propositions 2.2 statistic $\tilde{S}_n^h$ is bounded in probability under $H_0$, and by (2.2) $\kappa_n(h) \to \infty$ for every alternatives of $\Im_{alt}$. Hence $\alpha_n \to 0$, this convergence depends on the rate of convergence of $\varepsilon(n)$ to zero, that is it depends on the sub-family of considering alternatives, say $\Im \subseteq \Im_{alt}$.

We define IARE of $h$-test wrt $\psi$-test in harmony with idea of Inglot (2020). For each $n$ and every $k \geq 1$ let $u_k = u_k(n)$ be exact critical value of the $\psi$-test at the level $\alpha_n$ (2.5) and sample size $k$, i.e.,

$$P_0\{\tilde{S}_k^\psi > u_k\} = \alpha_n. \tag{2.6}$$

Next, for each $n$ let $k_n$ be the minimal sample size such that for all integer $j \geq 0$

$$P_1\{\tilde{S}_{k_n+j}^\psi > u_{k_n+j}\} \geq P_1\{\tilde{S}_n^h > z + \kappa_n(h)\}, \tag{2.7}$$

i.e. $k_n$ is the minimal sample size beginning from which the power of the $\psi$-test at the level $\alpha_n$ under alternatives of $\Im \subseteq \Im_{alt}$ is not smaller than that for the $h$-test at the same level and for the sample size $n$.

**Definition 2.1**. If there exists the limit

$$\lim_{n \to \infty} \frac{k_n}{n} = e \in [0, \infty],$$

then $e$ is called the IARE of $h$-test wrt $\psi$-test, for every alternatives of $\Im$. □

**Theorem 2.1.** Assume $h, \psi \in S$, the family of alternatives $\Im_{alt}$ is defined as in (2.1), and

(a) $N(\cdot)$ is taken as a positive function of the continuous variable $x$ is regularly varying of index $q \in (0, 2)$, viz., $N(ax) \sim a^q N(x)$ as $x \to \infty$, for all $a > 0$,

(b) $S_n^h$ satisfies condition (2.4),

(c) $\Im_1$, a sub-family of $\Im_{all}$, and a sequence of positive numbers $(\tau_n)$, $\tau_n \to 0$ or $\tau_n = c \leq 1/2$ for all $n$, such that for every sequence of positive numbers $(y_n)$, $y_n = O(\kappa_n(h))$, it holds

$$-\lim_{n \to \infty} \frac{1}{\tau_n y_n^2} \log P_0\{\tilde{S}_n^h > y_n\} = 1.$$

(d) There exists a sequence of positive numbers $(a_n)$, $a_n \to \infty$, such that

$$-\lim_{n \to \infty} \frac{2}{y_n^2} \log P_0\{\tilde{S}_n^\psi > y_n\} = 1,$$

for every sequence of positive numbers $(y_n)$ such that $y_n \to \infty$ and $y_n = o(a_n)$.

(f) $\Im_2$ is a sub-family of $\Im_{all}$ such that $\varepsilon(n) = o(a_n \sqrt{N}/n)$.

Then for every alternatives of $\Im = \Im_1 \cap \Im_2$ one has

$$\lim_{n\to\infty}\frac{k_n}{n}=\lim_{n\to\infty}\left(2\tau_n\frac{\rho^2(h,\lambda_n)}{\rho^2(\psi,\lambda_{k_n})}\right)^{1/(2-q)}<\infty. \tag{2.8}$$

□

**Remark 2.1**. Conditions (*c*) and (*d*) should be considered as conditions for the region, where it is supposed that under the hypothesis there are corresponding results on the probabilities of large deviations of the test statistics. In fact these regions define the family of alternatives $\Im_1$ and $\Im_2$, respectively.

### 2. 2. IARE of the chi-square test wrt *h*-tests.

Now on whenever it is convenient we shall use notation $d_n \ll b_n$ instead of well-known notation $d_n = o(b_n)$, also $\varepsilon(n) \ll d_n$ stands for the "sub-family of $\Im_{alt}$ is specified by the condition $\varepsilon(n) \ll d_n$".

We wish to evaluate the efficiency of the chi-square test by comparing its sensitivity relative to that of the sensitivity of $\psi$-test in terms of IARE, hence the chi-square test becomes as $h$-test in Theorem 2.1. Let $e(\chi_n^2, S_n^\psi,)$ denote the IARE of the chi-square test wrt to the $\psi$-test, $\psi_d$ be defined as in (1.3). We consider sparse, very sparse and dense models separately. Recall that we assume $N(an) \sim a^q N(n)$ as $n \to \infty$, where index $q \in (0,2)$.

*Sparse models*, that is $\lambda_n \to \lambda \in (0,\infty)$.

**Theorem 2.2**. Let $\lambda_n \to \lambda > 0$.

(i) If $\psi = \psi_d$, $d > -1$, $d \neq 1$ and $\varepsilon(n) \ll n^{-d^*/(1+2d^*)}$, where $d^* = \max(1,d)$, or $\psi$ satisfies the Cramèr condition $E\exp\{t|\psi(\xi)|\} < \infty$, for some $t > 0$, and $\varepsilon(n) \ll n^{-1/3}$, then

$$e(\chi_n^2, S_n^\psi) > 1.$$

(ii) If $\psi$ satisfies the Cramèr condition and $\varepsilon(n) \gg n^{-1/3}\log^{2/3} n$, Then

$$e(\chi_n^2, S_n^\psi) = 0. \qquad □$$

**Remark 2.2.** We emphasize that the Cramèr condition is fulfilled for CS (1.5), as well as for PDS $CR_n(d)$, i.e., for $\psi_d$, if $d \in (-1,0]$. For example, for the statistics $\Lambda_n$ and $T_n^2$, see (1.4), the Cramèr condition is fulfilled, but for $\chi_n^2$ it is not.

From Theorem 2.2 it follows

**Corollary 2.1**. Let $\lambda_n \to \lambda > 0$. Then for alternatives such that $\varepsilon(n) \ll n^{-1/3}$ it holds $e(\chi_n^2, \Lambda_n) > 1$, whereas for alternatives such that $\varepsilon(n) \gg n^{-1/3}\log^{2/3} n$ it holds $e(\chi_n^2, \Lambda_n) = 0$. □

Similar statement is valid for Freeman-Tukey statistic and each of CS (1.5).

*Very sparse models*, that is we assume that $\lambda_n \to 0$, but $n\lambda_n \to \infty$.

**Theorem 2.3**. Let $\lambda_n \to 0$, $n\lambda_n^3 \to \infty$.

(A) If $d \in (-1, 1/2)$ and $\lambda_n \gg n^{-(1-2d^*)/6}$, where now $d^* = \max(0, d)$, then

(i) for alternatives such that $\varepsilon(n) \ll n^{-1/3}$ it holds $e(\chi_n^2, S_n^{\psi_d}) = 1$,

(ii) for alternatives such that $(n\lambda_n)^{-1/3} \log^{2/3}(N^2/n) \ll \varepsilon(n) \ll (n\lambda_n^2)^{-1/4}$ it holds $e(\chi_n^2, S_n^{\psi_d}) = 0$.

(B) If $d \geq 1/2$ and alternatives such that $\varepsilon(n) \ll (\lambda_n^{1-d}/n^d)^{1/(2d+1)}$ then $e(\chi_n^2, S_n^{\psi_d}) = 1$. □

**Corollary 2.2**. Let $\lambda_n \to 0$, $n\lambda_n^3 \to \infty$ and $\lambda_n \gg n^{-1/6}$. Then for alternatives such that $\varepsilon(n) \ll n^{-1/3}$ it holds $e(\chi_n^2, \Lambda_n) = 1$, whereas for alternatives such that $(n\lambda_n)^{-1/3} \log^{2/3}(N^2/n) \ll \varepsilon(n) \ll (n\lambda_n^2)^{-1/4}$ it holds $e(\chi_n^2, \Lambda_n) = 0$. □

Actually, from Theorem 2.3 it follows alike to Corollary 2.2 result for every $\psi_d$-test, where $d \in (-1, 0]$. We emphasize that $n\lambda_n^3 \to \infty$ (i.e., $n \ll N \ll n^{4/3}$) is a necessary condition for the family of alternatives with $\varepsilon(n) \ll (\lambda_n^{1-d}/n^d)^{1/(2d+1)}$ to be a non-empty sub-family of $\Im_{alt}$.

**Theorem 2.4.** If $\lambda_n \to 0$ and $n\lambda_n \to \infty$. Then for alternatives such that $\varepsilon(n) \ll n^{-1/3}$ one has

$$e(\chi_n^2, \mu_r) = 1, \ r = 0, 1, 2, \ \text{and} \ \ e(\chi_n^2, C_n) = 1.$$ □

*Dense models* assume that $\lambda_n \to \infty$.

**Theorem 2.5.** Let $\lambda_n \to \infty$. (i) For every $d > -1$ and alternatives such that $\varepsilon(n) \ll (n\lambda_n^2)^{-1/3}$ it holds $e(\chi_n^2, S_n^{\psi_d}) = 1$. (ii) If $N = o(n^{3/8})$ then for every alternative of $\Im_{alt}$ it holds $e(\chi_n^2, \Lambda_n) = 1$.
□

The PDS can be rewritten as

$$CR_N(d) = \lambda_n \sum_{m=1}^{N} \bar{\psi}_d(\eta_m), \tag{2.9}$$

where

$$\bar{\psi}_d(x) = \frac{2}{d(d+1)}\left[(x/\lambda_n)^{d+1} - (d+1)(x/\lambda_n) + d\right], d \neq 0,$$

$$\bar{\psi}_0(x) = 2\left[(x/\lambda_n)\log(x/\lambda_n) - (x-\lambda_n)/\lambda_n\right].$$

Further, under $H_0$ the r.v. $(\eta_m - \lambda_n)/\sqrt{\lambda_n}$ has asymptotically normal distribution, as $\lambda_n \to \infty$, since $\eta_m \sim Bi(n, N^{-1})$. Hence $(\eta_m - \lambda_n)/\lambda_n = O_p(\lambda_n^{-1/2})$. Note that $\eta_m/\lambda_n = 1 + (\eta_m - \lambda_n)/\lambda_n$. Use these facts and Taylor expansion formula to get $\lambda_n \bar{\psi}_d(\eta_m) = (\eta_m - \lambda_n)^2/\lambda_n + O_p(\lambda_n^{-1/2}), d > -1$.

Thus, if $\lambda_n \to \infty$ then under $H_0$ each PDS asymptotically in distribution coincides with the chi-square statistic. This property is reflected to some extent in Theorem 2.5.

## 3. Proofs

We first present the proof of Proposition 2.1 and the theorems, after which we will give the proof of Lemma 2.1, 2.2 and 2.3,

**Proof of Proposition 2.1**. Recall $p_m = N^{-1}(1 + \varepsilon_m)$, $m = 1, 2, ..., N$, where $\varepsilon_1 + ... + \varepsilon_N = 0$, $\max_m |\varepsilon_m| \to 0$. Set $\xi \sim Poi(\lambda_n)$, $\varepsilon_j(n) = N^{-1}(\varepsilon_1^j + ... + \varepsilon_N^j)$, $\varepsilon_2(n) := \varepsilon(n)$,

$$\varphi_2(\xi) = \xi(\xi - 1) - 2\lambda_n \xi + \lambda_n^2 = (\xi - \lambda_n)^2 - (\xi - \lambda_n) - \lambda_n,$$

$$\varphi_3(\xi) = \xi(\xi - 1)(\xi - 2) - 3\lambda_n \xi(\xi - 1) + 3\lambda_n^2 \xi - \lambda_n^3$$
$$= (\xi - \lambda_n)^3 - 3(\xi - \lambda_n)^2 + (2 - 3\lambda_n)(\xi - \lambda_n) + 2\lambda_n.$$

Using Taylor expansion of $(1 + \varepsilon_m)^k \exp\{-\lambda_n \varepsilon_m\}$ we obtain

$$A_{1n}(h) = Eh(\xi) + \frac{1}{2}\varepsilon_2(n)E\big(h(\xi)\varphi_2(\xi)\big) + \frac{1}{6}\varepsilon_3(n)(E\big(h(\xi)\varphi_3(\xi)\big) + o(1)). \tag{3.1}$$

We have

$$E\varphi_2(\xi) = 0, \, E\xi\varphi_2(\xi) = 0, \, Var\varphi_2(\xi) = 2\lambda_n^2 \tag{3.2}$$

$$E\varphi_3(\xi) = 0, \, E\xi\varphi_3(\xi) = 0, \, Var\varphi_3(\xi) = 6\lambda_n^3.$$

Hence $E\big(h(\xi)\varphi_l(\xi)\big) = E\big(g(\xi)\varphi_l(\xi)\big)$. Due to these facts we obtain

$$A_{1n}(h) - A_{0n}(h) = \frac{\lambda_n \sigma(h)}{\sqrt{2}}\varepsilon_2(n)corr(g(\xi), \varphi_2(\xi))$$

$$+ \frac{\lambda_n^{3/2} \sigma(h)}{\sqrt{6}}\varepsilon_3(n)corr(g(\xi), \varphi_3(\xi))(1 + o(1)).$$

Therefore

$$\kappa_n(h) = \frac{\sqrt{n\lambda_n}}{\sqrt{2}}\varepsilon_2(n)corr(g(\xi), \varphi_2(\xi)) + \frac{\sqrt{n}}{\sqrt{6}}\lambda_n\varepsilon_3(n)corr(g(\xi), \varphi_3(\xi))(1 + o(1)).$$

Equality (2.2) follows, since $corr(g(\xi), \varphi_2(\xi)) = \rho(h, \lambda_n)$ and $\varepsilon_3(n) = o(\lambda_n^{-1/2})\varepsilon_2(n)$ by (2.1).

Next, by Taylor expansion

$$\sigma_{1n}^2(h) = Eg^2(\xi) + \frac{1}{2}\varepsilon_2(n)E\big(g^2(\xi)\varphi_2(\xi)\big)(1 + o(1)). \tag{3.3}$$

If $\lambda_n \to 0$ then $E\big(g^2(\xi)\varphi_2(\xi)\big) = \lambda_n^2(g^2(2) - 2g^2(1) + \sigma^2(h))(1 + o(1))$. If $\lambda_n \to \lambda \in (0, \infty)$ then evidently $E\big(g^2(\xi)\varphi_2(\xi)\big) \leq c(\lambda)$. By virtue of these facts and (3.3), equality (2.3) follows for very sparse and sparse models. Let $\lambda_n \to \infty$, then (3.3) gives

$$\sigma_{1n}^2(h) = \sigma^2(h)\left(1 + \frac{\lambda_n \varepsilon_2(n)}{2} E\left(\frac{g(\xi)}{\sigma(h)} \frac{\xi - \lambda_n}{\sqrt{\lambda_n}}\right)^2\right) = \sigma^2(h)(1 + o(1)),$$

since $\lambda_n \varepsilon_2(n) = o(1)$, $(\xi - \lambda_n)/\sqrt{\lambda_n} = O_p(1)$ and $Eg^2(\xi) = \sigma^2(h)$. Proof of Proposition 2.1 is completed. □

**Proof of Theorem 2.1.** By Proposition 2.2 we have

$$P_1\{\tilde{S}_n^h > z + \kappa_n(h)\} = P_1\{(\tilde{S}_{n,N}^h - NA_{1,n}(h))/\sqrt{N}\sigma_{1,n}(h) > z\sigma_{1,n}^{-1}(h)\sigma_{0,n}(h)\} = \Phi(-z) + o(1).$$

Hence, for arbitrary fixed $z$ and sequences of powers of $h$-test at the significance level defined by (2.5) we have

$$0 < c_1 \le P_1\{\tilde{S}_n^h > z + \kappa_n(h)\} \le c_2 < 1. \tag{3.4}$$

By virtue of this, definition of $k_n$, (2.3) and Proposition 2.2 we obtain

$$0 < c_1 \le P_1\{\tilde{S}_{k_n}^\psi > u_{k_n}\} = P_1\left\{\frac{S_{k_n}^\psi - NA_{1,k_n}(\psi)}{\sqrt{N}\sigma_{1,k_n}(\psi)} > (u_{k_n} - \kappa_{k_n}(\psi))(1 + o(1))\right\}. \tag{3.5}$$

Hence,

$$u_{k_n} - \kappa_{k_n}(\psi) \le c_3. \tag{3.6}$$

For every alternative of $\Im_{alt}$ we have $z + \kappa_n(h) \sim \kappa_n(h) \sim \sqrt{n\lambda_n/2}\,\varepsilon(n)\rho(h,\lambda_n) \to \infty$, and under alternatives $\Im_1$ by condition (c)

$$-\log \alpha_n = -\log P_0\{\tilde{S}_n^h > z + \kappa_n(h)\} \sim \frac{1}{2}\tau_n n \lambda_n \varepsilon^2(n) \rho^2(h,\lambda_n). \tag{3.7}$$

On the other hand $u_{k_n} \to \infty$, since (2.6) and $\alpha_n \to 0$, also $u_{k_n} = O(\kappa_{k_n}(\psi)) = o(a_{k_n})$ under alternatives $\Im_2$, by (3.6) and condition (d). Hence due to condition (d), choice of $u_{k_n}$, see (2.7), and (3.7) we obtain

$$-\log P_0\{\tilde{S}_{k_n}^\psi > u_{k_n}\} = \frac{1}{2}u_{k_n}^2(1 + o(1)) = -\log \alpha_n = \frac{1}{2}\tau_n n \lambda_n \rho^2(h,\lambda_n)\varepsilon^2(n)(1+o(1)),$$

i.e. $u_{k_n}^2 = \tau_n n \lambda_n \rho^2(h,\lambda_n)\varepsilon^2(n)(1+o(1))$. Thus, for every sequence of alternatives of $\Im$ we have

$$\tau_n n \lambda_n \rho^2(h,\lambda_n)\varepsilon^2(n) \le 2^{-1} k_n \lambda_{k_n} \rho^2(\psi,\lambda_{k_n})\varepsilon^2(k_n), \tag{3.8}$$

since (3.6) and definition of $\kappa_{k_n}(\psi)$.

Further, $P_1\{\tilde{S}_{k_n-1}^\psi > u_{k_n-1}\} \le P_1\{\tilde{S}_n^h > z + \kappa_n(h)\} \le c_2 < 1$ by definition of $k_n$ and (3.4). So

$$P_1\{\tilde{S}^{\psi}_{k_n-1} > u_{k_n-1}\} = P_1\left\{\frac{S^{\psi}_{k_n-1} - NA_{1,k_n-1}(\psi)}{\sqrt{N}\sigma_{1,k_n-1}(\psi)} > (u_{k_n-1} - \kappa_{k_n-1}(\psi))(1+o(1))\right\} \le c_2 < 1,$$

hence

$$u_{k_n-1} - \kappa_{k_n-1}(\psi) \ge -c_4. \tag{3.9}$$

So, $u_{k_n-1} \to \infty$, because $\kappa_{k_n-1}(\psi) \to \infty$ under alternatives $\Im_{all}$. Put

$$\kappa'_{k_n-1}(\psi) := \kappa_{k_n-1}(\psi) \frac{n\sqrt{N(k_n)}}{\rho(\psi, \lambda_{k_n-1})\sqrt{N(n)}(k_n-1)} = \frac{n}{\sqrt{2N}}\varepsilon(k_n).$$

Then $\kappa'_{k_n-1}(\psi) \to \infty$, $\kappa'_{k_n-1}(\psi) = o(a_n)$ under alternatives $\Im$. Hence by assumption (d) we obtain

$$\log P_0\{\tilde{S}^{\psi}_{k_n-1} > \kappa'_{k_n-1}(\psi)\} = -\frac{1}{2}\kappa'^2_{k_n-1}(\psi)(1+o(1)) = -\frac{1}{4}n\lambda_n\varepsilon^2(k_n)(1+o(1)). \tag{3.10}$$

On the other hand by (3.7) and definition of $u_{k_n-1}$ we have

$$\log P_0\{\tilde{S}^{\psi}_{k_n-1} > u_{k_n-1}\} = \log \alpha_n = -\frac{1}{2}\tau_n n\lambda_n \rho^2(h, \lambda_n)\varepsilon^2(n)(1+o(1))$$

$$\ge -\frac{1}{4}n\lambda_n \varepsilon^2(n)(1+o(1)) = \log P_0\{\tilde{S}^{\psi}_{k_n-1} > \kappa'_{k_n-1}\}, \tag{3.11}$$

since $\tau_n \le 1/2$, (3.10) and that it have to be $\varepsilon^2(n) \sim \varepsilon^2(k_n)$ for the alternatives $\Im$. From (3.11) it follows that $u_{k_n-1} \le \kappa'_{k_n-1} = o(a_n)$. Then by assumption (d)

$-\log P_0\{\tilde{S}^{\psi}_{k_n-1} > u_{k_n-1}\} = 2^{-1}u^2_{k_n-1}(1+o(1))$. Due to this fact, (3.9), (3.10) and (3.11) we obtain

$$2^{-1}(k_n-1)\lambda_{k_n-1}\rho^2(\psi, \lambda_{k_n})\varepsilon^2(k_n)(1+o(1)) = 2^{-1}\kappa^2_{k_n-1}(\psi) \le 2^{-1}u^2_{k_n-1}(1+o(1))$$

$$= -\log P_0\{\tilde{S}^{\psi}_{k_n-1} > u_{k_n-1}\} = -\log \alpha_n = \tau_n n\lambda_n \rho^2(h, \lambda_n)\varepsilon^2(n)(1+o(1)).$$

From this, inequality (3.8) and that $\varepsilon^2(n) \sim \varepsilon^2(k_n)$, we have for the alternatives $\Im$

$$2^{-1}(k_n-1)\lambda_{k_n-1}\rho^2(\psi, \lambda_{k_n})(1+o(1)) \le \tau_n n\lambda_n \rho^2(h, \lambda_n) \le 2^{-1}k_n\lambda_{k_n}\rho^2(\psi, \lambda_{k_n}).$$

That is reminding that $\lambda_l = l/N(l)$ we obtain

$$\frac{N(n)}{n^2} \sim \frac{N(k_n)}{k_n^2}b_n, \text{ where } b_n = 2\tau_n \frac{\rho^2(h, \lambda_n)}{\rho^2(\psi, \lambda_{k_n})}. \tag{3.12}$$

Note that $\lim_{n\to\infty} b_n = b < \infty$, because $\psi \in S$. From (3.12) due to condition (a) we have

$$\frac{N(n)}{n^2} \sim \frac{N(a_n k_n)}{(a_n k_n)^2}, \text{ where } a_n = b_n^{-1/(2-q)}. \tag{3.13}$$

Now we use arguments similar to those of Quine and Robinson (1985, page 730). The function $R(x) = N(x)/x^2$ is regularly varying with index $q-2 < 0$, therefore if $0 < r_x < \infty$ and $r_x$ has a limit in $[0,\infty]$ as $x \to \infty$ then $R(r_x x)/R(x) \sim r_x^{q-2}$ as $x \to \infty$. Let a sub-sequence of $k_n/n$ have a further sub-sequence for which $k_n/n \to d \in [0,\infty]$, then for this subsequence $R(a_n k_n)/R(n) = R((a_n k_k/n)n)/R(n) \sim (a_n d)^{q-2}$. On the other hand $R(a_n k_n)/R(n) \sim 1$ since (3.13), and hence $(a_n d)^{q-2} \sim 1$, that is $d = b^{1/(2-q)}$ do not depend on particular sub-sequence. Thus for whole sequence $k_n/n \to b^{1/(2-q)}$. Theorem 2.1 follows. □

**Proof of Theorems 2.2 -2.5.** A few remarks first. We use the chi-square test as a benchmark and it is considered as $h$-test of Theorem 2.1, hence the condition (b) is fulfilled. Also, since $\rho(h, \lambda_n) = 1$ if $h(x) = x^2$ and $\psi \in S$, that is $|\rho(\psi, \lambda_n)| \geq c > 0$, it follows from Theorem 2.1 that under the conditions (a),(c),(d) and (f): if $\tau_n = c$ then there exists a constant $c_1 > 0$ such that $\lambda_{k_n} \sim \lambda_n(k_n/n)(N(n)/N(k_n)) \sim \lambda_n(k_n/n)^{1-q} \sim c_1 \lambda_n$, also if $\tau_n \to 0$ then $e(\chi_n^2, S_n^\psi) = 0$. We will use these facts. Finally, it suffices to indicate that conditions (c), (d) and (f) are satisfied under the conditions of each of Theorems 2.2–2.5, so we will no longer refer to Theorem 2.1 in the course of proving the theorems. All relevant large deviations results for the chi-square statistic are collected in the following Proposition 3.1. Note that for the chi-square statistic

$$h(u) = (u - \lambda_n)^2/\lambda_n, \ Eh(\xi) = 1, \ \sigma^2(h) = 2, \ \rho(h, \lambda_n) = 1, \ \kappa_n(h) = \kappa_n = \sqrt{n\lambda_n/2}\varepsilon(n).$$

**Proposition 3.1**. (A) Let either

(i) $\lambda_n \to \lambda \in (0,\infty)$ and $\varepsilon(n) \ll n^{-1/3}$, or

(ii) $\lambda_n \to 0$, $n\lambda_n^3 \to \infty$ and $\varepsilon(n) \ll n^{-1/3}$, or

(iii) $\lambda_n \to \infty$ and $\varepsilon(n) \ll (n\lambda_n^2)^{-1/3}$, or

(iv) $N \ll n^{3/8}$ and $\sqrt{n\lambda_n}\varepsilon(n) \to \infty$ and $\varepsilon(n) \ll \lambda_n^{-1}$.

Then for every sequence of $y_n = O(\sqrt{n\lambda_n}\varepsilon(n))$ it holds

$$-\lim_{n \to \infty} \frac{2}{y_n^2} \log P_0\left\{(\chi_N^2 - N)/\sqrt{2N} > y_n\right\} = 1.$$

(B) Let $\sqrt{n} \leq N$ and

$$\varepsilon(n) = (n\lambda_n)^{-1/3}\omega_n^{2/3}, \text{ where } \omega_n \to \infty, \omega_n = o(\sqrt{n\lambda_n}), \log(N^2/n) = o(\omega_n). \tag{3.14}$$

Then there exists a sequence $(\tau_n)$, $\tau_n \to 0$ such that for every sequence $u_n$ such that $u_n = O(\sqrt{n\lambda_n}\varepsilon(n))$ it holds

$$-\lim_{n \to \infty} \frac{1}{\tau_n u_n^2} \log P_0\left\{(\chi_n^2 - N)/\sqrt{2N} > u_n\right\} = 1. \qquad \square$$

**Remark 3.1**. In fact, there is no need for $\varepsilon(n)$ to be the same as in (1.1). Proposition 3.1 still true for an arbitrary sequence of positive numbers $\varepsilon(n) \to 0$ that satisfies the corresponding conditions.

**Proof of Proposition 3.1**. By Corollary 4.2 of Mirakhmedov (2022) it holds

$$\log P_0 \left\{ \chi_N^2 > x_n \sqrt{2N} + N \right\} = -\frac{1}{2} x_n^2 (1+o(1)) \qquad (3.15)$$

for arbitrary $\lambda_n$ and every sequence of positive numbers $x_n \to \infty$ such that

$$x_n = o\left( (\sqrt{N} \min(1, \lambda_n^2))^{1/3} \right). \qquad (3.16)$$

Also, Corollary 2.6, Eq (2.17), of Kallenberg (1985) gives that (3.15) still hold if

$$N \ll \sqrt{n}, \ x_n = o(\sqrt{N}) \text{ and } N^{-3/2} n^{1/2} x_n \to \infty. \qquad (3.17)$$

Part (A) of Proposition 3.1 follows if we choose $\varepsilon(n)$ so that $\kappa_n = n\varepsilon(n)/\sqrt{2N}$ satisfies (3.16) and (3.17) in each cases, respectively. Indeed, for the sparse and dense models (cases (i) and (iii)) we have in (3.16) $x_n \ll N^{1/6}$, whereas for very sparse models (case (ii)) $x_n \ll (\sqrt{N}\lambda_n^2)^{1/3} = (n\lambda_n^3)^{1/6}$. Using these facts we can easily observe that if $\varepsilon(n)$ is as in the assertions (i),(ii) and (iii) then for every positive $y_n = O(\kappa_n)$ Eq. (3.16) is valid. Part (iv) assumes that we are dealing with dense model, and hence we can restrict ourselves by $\varepsilon(n) \geq (n\lambda_n^2)^{-1/3}$, since case (iii). Then (3.17) is fulfilled with $x_n = \kappa_n$ if $N \ll n^{3/8}$. Part (A) is proven completely. Proof of Part (B) is alike to proof of part (iii) of Theorem 3.4 of Mirakhmedov (2021, p.16-18). □

**Proposition 3.2**. Let $\lambda_n \to \lambda \in (0, \infty)$ and $d^* = \max(1, d)$ For every $d > -1$ and $x_n \to \infty$, $x_n = o(n^{1/2(1+2d^*)})$ one has

$$-\log P_0 \left\{ \tilde{S}_n^{\psi_d} > x_n \right\} = \frac{x_n^2}{2}(1+o(1)). \qquad (3.18)$$

Moreover, (3.18) remains true for $x_n \to \infty$, $x_n = o(\sqrt{n})$, if statistics $S_n^\psi$ satisfy the Cramèr condition.

**Proof.** By Corollary 3.1 of Mirakhmedov (2022) Eq.(3.18) is valid if $d > 0$ and $x_n = o(\min(n^{1/6}, n^{1/(1+2d)})) = o(n^{1/(1+2d^*)})$. Further, it follows from Theorem 2 of Ivchenko and Mirakhmedov (1995) that Eq.(3.18) holds for $x_n \to \infty$, $x_n = o(\sqrt{n})$ if statistics $S_n^\psi$ satisfy the Cramèr condition, in particular for PDS with tuning parameter $d \in (-1, 0]$. Proposition 3.2 follows.

□

**Proof of Theorem 2.2**. Condition (*c*), where $\mathfrak{I}_1$ such that $\varepsilon(n) \ll n^{-1/3}$ and $\tau_n = 1/2$, is fulfilled since Proposition 3.1 part (A) case (i). Condition (*d*) is fulfilled with $a_n = n^{1/2(1+2d^*)}$ by Proposition 3.2, hence the family $\mathfrak{I}_2$ of condition (f) we should determine so that $n\varepsilon(n)/\sqrt{2N} \ll n^{1/2(1+2d^*)}$, that is $\varepsilon(n) \ll n^{-d^*/(1+2d^*)}$. Then $\mathfrak{I} = \mathfrak{I}_1 \cap \mathfrak{I}_2$ is such that $\varepsilon(n) \ll n^{-d^*/(1+2d^*)}$. Part (i) follows, since in this case $|\rho(\psi, \lambda_{k_n})| < 1$. Further, under the conditions of part (ii) $\mathfrak{I}_1$ is such that $\varepsilon(n) \gg n^{-1/3} \log^{2/3} n$ and $\tau_n \to 0$ by Proposition 3.1 part (B), whereas by Proposition 3.2 in this case $a_n = \sqrt{n}$, and hence $\mathfrak{I}_2 = \mathfrak{I}_{alt}$. Part (ii) follows. □

**Proof of Theorem 2.3.** Set $W_N(d) = \min(n^{1/4}, (n\lambda_n^3)^{1/2(1+2d^*)})$, where $d^* = \max(0, d)$.

**Proposition 3.3**. Let $\lambda_n \to 0$ and $n\lambda_n^3 \to \infty$. Then uniformly in $x_n \geq 0$, $x_n = o(W_N(d))$ it holds

$$-\log P_0\{\tilde{S}_n^{\psi_d} > x_n\} = \frac{x_n^2}{2}(1+o(1)), \ d > -1.$$

**Proof** follows from Corollary 3.3 of Mirakhmedov (2022), because standardized r.v. $\tilde{S}_n^{h_d}$ coincides with accordingly standardized version of the r.v. $R_n^d = \eta_1^{d+1} + ... + \eta_N^{d+1}$, $d \neq 0$, and $R_n^0 = 2\eta_1 \ln \eta_1 + ... + 2\eta_N \ln \eta_N$.

*Proof of part* (A). By Proposition 3.1 (ii) Condition (c) is fulfilled, where $\mathfrak{I}_1$ defines as $\varepsilon(n) \ll n^{-1/3}$ and $\tau_n = 1/2$. Next, by Proposition 3.3 we have $a_n = n^{1/4}$, because if $N \leq n^{(7-2d^*)/6}$ and $d \in (-1, 1/2)$ then $W_n(d) = n^{1/4}$. So, $\mathfrak{I}_2$ defines by $\varepsilon(n) \ll \sqrt{N}/n^{3/4} = (n\lambda_n^2)^{-1/4}$. Then $\mathfrak{I} = \mathfrak{I}_1$, since $n \ll N$. Part (i) follows, because $\rho(\psi_d, \lambda_n) \to 1$ by Lemma 2.2. Part (ii) also follows by similar reasons using Proposition 3.1 part (B) (instead of part (ii)) due to which $\tau_n \to 0$ and $\mathfrak{I}_1$ is defined by $\varepsilon(n) = (n\lambda_n)^{-1/3} \omega_n^{2/3} \ll (n\lambda_n^2)^{-1/4}$, where $\omega_n \ll (N/\sqrt{n})^{1/4}$ (note that $(N/\sqrt{n})^{1/4} \leq n/\sqrt{N}$, because $N \ll n^{4/3}$ by condition $n\lambda_n^3 \to \infty$).

*Proof of part* (B). From Proposition 3.1 and 3.3 we have: $\mathfrak{I}_1$ is such that $\varepsilon(n) \ll n^{-1/3}$, $a_n = W_N(d) = (n\lambda_n^3)^{1/2(2d+1)}$ since $d \geq 1/2$, and $\mathfrak{I}_2$ is such that $\varepsilon(n) \ll [n^d \lambda_n^{d-1}]^{-1/(2d+1)}$, $\tau_n = 1/2$. Hence $\mathfrak{I}_2 \subset \mathfrak{I}_1$ because $N \ll n^{4/3}$. Part (B) follows since Lemma 2.2. □

**Proof of Theorem 2.4.** Set $\pi_r(\lambda) = \lambda^r e^{-\lambda}/r!$,

$$A_{rn} = \sum_{m=1}^{N} \pi_r(\lambda_n), \quad \sigma_{rn}^2 = \sum_{m=1}^{N} \pi_r(\lambda_n)(1-\pi_r(\lambda_n)) - n^{-1}\left(\sum_{m=1}^{N}(r-\lambda_n)\pi_r(\lambda_n)\right)^2.$$

From Corollary 4.3 of Mirakhmedov (2022) it follows

**Proposition 3.4.** Let $\lambda_n \to 0$ and $n\lambda_n \to \infty$. Then for $r = 0, 1, 2$ and $x_n \to \infty$, $x_n = o(n^{1/4})$ it holds

$$\log P_0\{\mu_r > \sigma_{r,n} x_n + A_{r,n}\} = -x_n^2/2(1+o(1)).$$

By Proposition 3.1 and 3.4 we have: $\Im_1$ is such that $\varepsilon(n) \ll n^{-1/3}$, whereas $a_n = n^{1/4}$, $\tau_n = 1/2$, $\Im_2$ defines by $\varepsilon(n) \ll \sqrt{N}/n^{3/4} = (n\lambda_n^2)^{-1/4}$, and hence $\Im_1 \subset \Im_2$. Proof concludes, since Lemma 2.2, according of which $\rho(h, \lambda_n) \to 1$ for $h(x) = I\{x = r\}$ if $r = 0, 1, 2$. The second statement follows because $C_n = \mu_0 - (n - N)$. □

**Proof of Theorem 2.5.** We will use the following propositions.

**Proposition 3.5.** For every $d > -1$ and $x_n \to \infty$, $x_n = o(N^{1/6})$ one has

$$-\log P_0\{\tilde{S}_n^{\psi_d} > x_n\} = x_n^2/2(1 + o(1)).$$

**Proof** follows from Corollary 3.2 of Mirakhmedov (2022).

**Proposition 3.6.** If $N \ll n^{3/7}$ then for $x_n \to \infty$ and $x_n = o(\sqrt{N})$ it holds

$$-\log P_0\{\Lambda_n > x_n\sqrt{2N} + N\} = x_n^2/2(1 + o(1)).$$

**Proof** follows from Proposition 3.5 if $x_n = o(N^{1/6})$. If $x_n \geq N^{1/6}$, then by Eq. (2.13) of Kallenberg (1985) and that $N = o(n^{3/7})$ we have

$$\log P_0\{\Lambda_n > x_n\sqrt{2N} + N\} = -\frac{x_n^2}{2} + O\left(\frac{x_n^3}{\sqrt{N}} + \log N + \frac{N^{3/2}}{\sqrt{n}}\right) = -\frac{x_n^2}{2}(1 + o(1)).$$

To prove part (i) we define for every $d > -1$ the family $\Im_1$ by $\varepsilon(n) \ll (n\lambda_n^2)^{-1/3}$ according to Proposition 3.1(iii), and we set $\Im_2 = \Im_1$ because due to Proposition 3.5 $a_n = N^{1/3}$ and hence $\varepsilon(n) \ll N^{1/6}\sqrt{N}/n = (n\lambda^2)^{-1/3}$. To prove part (ii) we use Proposition 3.1 part (iv) and Proposition 3.6, according to which we define $\Im_1 = \Im_2 = \Im_{alt}$, because $a_n = \sqrt{N}$. Theorem 2.5 follows, since $\lambda_{k_n} \to \infty$ and hence $\rho(\psi_d, \lambda_{k_n}) \to 1$ by Lemma 2.3. □

**Proof of Lemma 2.1.** Consider two statistics $S_n^{h_1}$ and $S_n^{h_2}$. Denote

$$\tilde{\xi}_m = (\xi_m - \lambda_n)/\sqrt{n}, \ \tau_{jn} = \lambda_n^{-1}\mathrm{cov}(h_j(\xi), \xi), \ g_j(x) = h_j(x) - Eh_j(\xi) - \tau_{jn}(x - \lambda_n),$$

$$R_n^{g_j} = \sum_{m=1}^{N} g_j(\eta_m), \ V_n^{g_j} = \sum_{m=1}^{N} g_j(\xi_m), \ \zeta_n = \sum_{m=1}^{N}(\xi_m - \lambda_n).$$

We have $S_n^{h_j} = R_n^{g_j} + NEg_j(\xi)$ and $\mathrm{cov}(S_n^{h_1}, S_n^{h_2}) = \mathrm{cov}(R_n^{g_1}, R_n^{g_2})$. The following equality holds:

$$\mathrm{cov}_0(S_n^{h_1}, S_n^{h_2}) = \upsilon_n \int_{-\pi\sqrt{n}}^{\pi\sqrt{n}} E_0\left(V_n^{g_1} V_n^{g_1} \exp\{i\tau\zeta_n/\sqrt{n}\}\right) d\tau, \quad (3.19)$$

where

$$\upsilon_n := \left(2\pi\sqrt{n}P\{\zeta_n = 0\}\right)^{-1} = \frac{n!e^n}{2\pi n^n \sqrt{n}} = \frac{1}{\sqrt{2\pi}}\left(1 + o\left(\frac{1}{n}\right)\right).$$

**Indeed**. Note that $E_0 V_n^{g_l} = 0$. It is well known that $\mathcal{L}((\eta_1,...,\eta_N)) = \mathcal{L}((\xi_1,...,\xi_N)/\zeta_n = 0)$, where $\mathcal{L}(X)$ stands for the distribution of a random vector $X$. Hence $\operatorname{cov}_0(R_n^{g_1}, R_n^{g_2}) = E_0(V_n^{g_1} V_n^{g_2} | \zeta_n = 0)$. On the other hand $E(V_n^{g_1} V_n^{g_2} e^{i\tau\zeta_n}) = E\{e^{i\tau\zeta_n} E(V_n^{g_1} V_n^{g_2} | \zeta_h)\}$. Now (3.19) follows by Fourier inversion.

We have

$$\int_{-\pi\sqrt{n}}^{\pi\sqrt{n}} E_0\left(V_n^{g_1} V_n^{g_2} \exp\{i\tau\zeta_n/\sqrt{n}\}\right) d\tau = \sum_{m=1}^{N}\sum_{j=1}^{N} \int_{-\pi\sqrt{n}}^{\pi\sqrt{n}} E_0\left(g_1(\xi_m) g_2(\xi_j) \exp\{i\tau\zeta_n/\sqrt{n}\}\right) d\tau$$

$$= \sum_{m=1}^{N}\sum_{j=1}^{N} \int_{-\pi\sqrt{n}}^{\pi\sqrt{n}} E_0\left(g_1(\xi_m) g_2(\xi_j) \exp\{i\tau(\tilde{\xi}_m + \tilde{\xi}_j)\}\right) \prod_{l=1, l\neq m, j}^{N} E_0 \exp\{i\tau\tilde{\xi}_l\} d\tau$$

$$= \sum_{m=1}^{N} \int_{-\pi\sqrt{n}}^{\pi\sqrt{n}} \left(E_0 \exp\{i\tau\tilde{\xi}_1\}\right)^{N-1} E_0\left(g_1(\xi_m) g_2(\xi_m) e^{i\tau\tilde{\xi}_m}\right) d\tau$$

$$+ \sum_{\substack{m,j=1\\m\neq j}}^{N} \int_{-\pi\sqrt{n}}^{\pi\sqrt{n}} \left(E_0 \exp\{i\tau\tilde{\xi}_1\}\right)^{N-2} E_0 g_1(\xi_m) e^{i\tau\tilde{\xi}_m} E_0 g_2(\xi_j) e^{i\tau\tilde{\xi}_j} d\tau$$

$$= N \int_{-\pi\sqrt{n}}^{\pi\sqrt{n}} \left(E_0 \exp\{i\tau\tilde{\xi}_1\}\right)^{N-1} E_0\left(g_1(\xi_1) g_2(\xi_1) e^{i\tau\tilde{\xi}_1}\right) d\tau$$

$$+ N(N-1) \int_{-\pi\sqrt{n}}^{\pi\sqrt{n}} \left(E_0 \exp\{i\tau\tilde{\xi}_1\}\right)^{N-2} E_0 g_1(\xi_1) e^{i\tau\tilde{\xi}} E_0 g_2(\xi_1) e^{i\tau\tilde{\xi}_1} d\tau =: NQ_1 + N(N-1)Q_2. \qquad (3.20)$$

We have

$$E_0 e^{i\tau\tilde{\xi}_1} = \exp\left\{\lambda_n\left(e^{i\tau/\sqrt{n}} - 1 - i\tau/\sqrt{n}\right)\right\} = \exp\left\{-\frac{\tau^2}{2N} + \frac{\theta|\tau|^3}{6N\sqrt{n}}\right\}, \qquad (3.21)$$

$$\left|E_0 e^{i\tau\tilde{\xi}_1}\right| = \exp\left\{-2\lambda_n \sin^2\frac{\tau}{2\sqrt{n}}\right\} \leq \exp\left\{-\frac{2\tau^2}{N\pi^2}\right\}. \qquad (3.22)$$

Write

$$Q_1 = \int_{-\pi\sqrt{n}}^{\pi\sqrt{n}} \exp\left\{-\frac{\tau^2(N-1)}{2N}\right\} E_0\left(g_1(\xi_1) g_2(\xi_1) e^{i\tau\tilde{\xi}_1}\right) d\tau$$

$$+ \int_{-\pi\sqrt{n}}^{\pi\sqrt{n}} \left[\left(E_0 \exp\{i\tau\tilde{\xi}_1\}\right)^{N-1} - \exp\left\{-\frac{\tau^2(N-1)}{2N}\right\}\right] E_0\left(g_1(\xi_1) g_2(\xi_1) e^{i\tau\tilde{\xi}_1}\right) d\tau =: Q_{11} + Q_{12}. \qquad (3.23)$$

We have

$$Q_{11} = E_0\left(g_1(\xi_1) g_2(\xi_1)\right) \int_{-\pi\sqrt{n}}^{\pi\sqrt{n}} \exp\left\{-\frac{\tau^2(N-1)}{2N}\right\} d\tau$$

$$+ \int_{-\pi\sqrt{n}}^{\pi\sqrt{n}} \exp\left\{-\frac{\tau^2(N-1)}{2N}\right\} E_0\left(g_1(\xi_1)g_2(\xi_1)(e^{i\tau\tilde{\xi}_1}-1)\right)d\tau$$

$$= E_0\left(g_1(\xi_1)g_2(\xi_1)\right)\sqrt{2\pi}(1+o(1))+cE_0\left|g_1(\xi_1)g_2(\xi_1)\tilde{\xi}_1^2\right|. \tag{3.24}$$

Also,

$$Q_{12} = \frac{c}{\sqrt{n}} E_0\left|g_1(\xi_1)g_2(\xi_1)\right| \int_{|\tau|\leq\sqrt{n}} |\tau|^3 \exp\left\{-\frac{\tau^2(N-1)}{3N}\right\}d\tau$$

$$+ E_0\left|g_1(\xi_1)g_2(\xi_1)\right| \int_{\sqrt{n}\leq|\tau|\leq\pi\sqrt{n}} \left(\exp\left\{-\frac{2\tau^2(N-1)}{\pi^2 N}\right\} + \exp\left\{-\frac{\tau^2(N-1)}{2N}\right\}\right)d\tau$$

$$= \frac{c}{\sqrt{n}} E_0\left|g_1(\xi)g_2(\xi_1)\right|(1+o(1)). \tag{3.25}$$

Use the facts that $E_0 g_i(\xi_1) = 0$, $E_0 g_i(\xi_1)\tilde{\xi}_1 = 0$ and inequality $\left|e^{it}-1-it\right|\leq t^2/2$ to get

$$|Q_2| = \left|\int_{-\pi\sqrt{n}}^{\pi\sqrt{n}} \left(E\exp\{i\tau\tilde{\xi}_1\}\right)^{N-2} E_0 g_1(\xi_1)(e^{i\tau\tilde{\xi}}-1-i\tau\tilde{\xi}_1)E_0 g_2(\xi_1)(e^{i\tau\tilde{\xi}_1}-1-i\tau\tilde{\xi}_1)d\tau\right|$$

$$\leq \frac{1}{4}\int_{-\pi\sqrt{n}}^{\pi\sqrt{n}} \left|E_0\exp\{i\tau\tilde{\xi}_1\}\right|^{N-2} \tau^4 E_0(|g_1(\xi_1)|\tilde{\xi}_1^2)E_0(|g_2(\xi_1)|\tilde{\xi}_1^2)d\tau$$

$$\leq \frac{1}{4}E_0(|g_1(\xi_1)|\tilde{\xi}_1^2)E_0(|g_2(\xi_1)|\tilde{\xi}_1^2)\int_{-\pi\sqrt{n}}^{\pi\sqrt{n}} \tau^4 \exp\left\{-\frac{2(N-2)\tau^2}{N\pi^2}\right\}d\tau$$

$$\leq cE_0(|g_1(\xi_1)|\tilde{\xi}_1^2)E_0(|g_2(\xi_1)|\tilde{\xi}_1^2) \leq c\left(E_0 g_1^2(\xi_1)\right)^{1/2}\left(E_0 g_2^2(\xi_1)\right)^{1/2} E_0\tilde{\xi}_1^4$$

$$= c\left(E_0 g_1^2(\xi_1)\right)^{1/2}\left(E_0 g_2^2(\xi_1)\right)^{1/2}(N^{-2}+3(nN)^{-1}). \tag{3.26}$$

Thus, by (3.19) - (3.26) we have

$$\frac{1}{N}\text{cov}_0\left(R_n^{g_1},R_n^{g_2}\right) = \text{cov}_0\left(g_1(\xi_1),g_2(\xi_1)\right)(1+o(1))$$

$$+cE_0\left|g_1(\xi_1)g_2(\xi_1)\tilde{\xi}_1^2\right| + \frac{c}{\sqrt{n}}E_0\left|g_1(\xi_1)g_2(\xi_1)\right|(1+o(1))+c\left(E_0 g_1^2(\xi_1)\right)^{1/2}\left(E_0 g_2^2(\xi_1)\right)^{1/2}(N^{-1}+3n^{-1}).$$

Use here the facts that $N^{-1}\text{Var}_0 R_n^{g_l} = E g_l^2(\xi)(1+o(1))$, $E|g_1(\xi)g_2(\xi)|\leq\sqrt{E g_1^2(\xi)E g_2^2(\xi)}$,

$\tilde{\xi}_1^2 = (\xi_1-\lambda_n)^2/N\lambda_n = o_p(N^{-1/2})$ to get

$$\text{corr}_0\left(R_n^{g_1},R_n^{g_2}\right) = \frac{N^{-1}\text{cov}_0\left(R_n^{g_1},R_n^{g_2}\right)}{\sqrt{N^{-1}\text{Var}_0 R_n^{g_1}\cdot N^{-1}\text{Var}_0 R_n^{g_2}}}$$

$$= \text{corr}\left(g_1(\xi),g_2(\xi)\right)(1+o(1))+O\left(\frac{1}{\sqrt{N}}+\frac{1}{\sqrt{n}}\right).$$

Lemma 2.1 follows by taking here $h_1 = h$ and $S_n^{h_2} = \chi_n^2$. □

**Proof of Lemma 2.2.** Set $\Delta h(x) = h(x+1) - h(x)$. Since $\lambda_n \to 0$ we have

$$Eh(\xi) = \left(1 - \lambda_n + \frac{1}{2}\lambda_n^2 + O(\lambda_n^3)\right)\left(h(0) + \lambda h(1) + \frac{1}{2}\lambda^2 h(2) + O(\lambda_n^3)\right)$$

$$= h(0) + \lambda_n \Delta h(0) + \frac{1}{2}\lambda_n^2 \Delta^2 h(0) + O(\lambda_n^3), \qquad (3.27)$$

$$Eh^2(\xi) = h^2(0) + \lambda \Delta h^2(0) + \frac{1}{2}\lambda^2 \Delta^2 h^2(0) + O(\lambda_n^3), \qquad (3.28)$$

$$(Eh(\xi))^2 = h^2(0) + 2\lambda h(0)\Delta h(0) + \lambda^2 (\Delta h(0))^2 + \lambda^2 h(0)\Delta^2 h(0) + O(\lambda_n^3), \qquad (3.29)$$

Further,

$$Eh(\xi)\xi = \left(1 - \lambda_n + \frac{1}{2}\lambda_n^2 + O(\lambda_n^3)\right)\left(\lambda_n h(1) + \lambda_n^2 h(2) + \lambda_n^3 2^{-1} h(3) + O(\lambda_n^4)\right)$$

$$= \lambda_n h(1) + \lambda_n^2 \Delta h(1)) + \lambda_n^3 2^{-1} \Delta^2 h(1) + O(\lambda_n^4). \qquad (3.30)$$

Use this and (3.27) to get

$$\tau_n = \lambda^{-1} \operatorname{cov}(h(\xi), \xi) = \lambda^{-1} Eh(\xi)\xi - Eh(\xi) = \Delta h(0) + \lambda \Delta^2 h(0).$$

Hence

$$\lambda \tau_n^2 = \lambda (\Delta h(0))^2 + 2\lambda^2 \Delta h(0)\Delta^2 h(0) + \lambda_n^3 [(\Delta^2 h(0))^2 + \Delta h(0)\Delta^3 h(0)]. \qquad (3.31)$$

By (3.28),(3.29) and (3.31) we obtain

$$\sigma^2(h) = Eh^2(\xi) - (Eh(\xi))^2 - \lambda \tau_n^2 = \lambda\left(\Delta h^2(0) - 2h(0)\Delta h(0) - (\Delta h(0))^2\right)$$

$$+ \frac{1}{2}\lambda^2\left(\Delta^2 h^2(0) - 2(\Delta h(0))^2 - 2h(0)\Delta^2 h(0) - 4\Delta h(0)\Delta^2 h(0)\right)$$

$$+ \frac{1}{6}\lambda^3\left[\Delta^3 h^2(0) - 6\Delta h(0)\Delta^2 h(0) - 2h(0)\Delta^3 h(0) - 6(\Delta^2 h(0))^2 - 6\Delta h(0)\Delta^3 h(0)\right]$$

$$:= \lambda_n \nabla_1 + \frac{1}{2}\lambda_n^2 \nabla_2 + \frac{1}{6}\lambda_n^3 \nabla_3. \qquad (3.32)$$

We have

$$\nabla_1 = \left(h^2(1) - h^2(0) - 2h(0)h(1) + 2h^2(0) - h^2(1) + 2h(1)h(0) - h^2(0)\right) = 0,$$

$$\nabla_2 = h^2(2) + 4h^2(1) + h^2(0) - 4h(0)h(1) + 2h(0)h(2) - 4h(1)h(2) = \frac{1}{2}\lambda^2 (\Delta^2 h(0))^2,$$

$$\nabla_3 = \Delta^3 h^2(0) - 6\Delta h(1)\Delta^2 h(0) - 2\Delta^3 h(0)[h(1) + 2\Delta h(0)].$$

Apply these in (3.32) to get

$$\sigma^2(h) = \frac{1}{2}\lambda^2 (\Delta^2 h(0))^2 \left[1 + \frac{\lambda_n}{3}\frac{\Delta^3 h^2(0) - 6\Delta h(1)\Delta^2 h(0) - 2\Delta^3 h(0)[h(1) + 2\Delta h(0)]}{(\Delta^2 h(0))^2} + o(\lambda_n)\right]$$

$$:= \frac{1}{2}\lambda^2(\Delta^2 h(0))^2[1+\lambda_n 3^{-1}\omega(h)] \quad \text{(for short)}. \tag{3.33}$$

We find that

$$Eh(\xi)\xi^2 = \lambda h(1) + \lambda^2[h(2)+\Delta h(1)] + \frac{1}{2}\lambda^3[2\Delta h(2)+\Delta^2 h(1)] + O(\lambda_n^4).$$

Using this, (3.30) and (3.27) after some algebra we obtain

$$\text{cov}\left(h(\xi)-\tau_n\xi, \xi^2-(2\lambda_n+1)\xi\right) = Eh(\xi)\xi^2 - (2\lambda_n+1)Eh(\xi)\xi + \lambda^2 Eh(\xi)$$

$$= \lambda_n^2 \Delta^2 h(0) + \frac{1}{2}\lambda_n^3[\Delta^3 h(1)+\Delta h(2)] + O(\lambda_n^4).$$

Due to this, (3.33) and (3.2) we obtain

$$\rho(h,\lambda_n) = \frac{\lambda_n^2\Delta^2 h(0) + \lambda_n^3 2^{-1}[\Delta^3 h(1)+\Delta h(2)]}{\lambda^2\Delta^2 h(0)\sqrt{1+\lambda_n 3^{-1}\omega(h)}}$$

$$= \left[\Delta^2 h(0)\right]^{-1}\left(\Delta^2 h(0)+\lambda_n 2^{-1}[\Delta^3 h(1)+\Delta h(2)]\right)\left(1-\lambda_n 6^{-1}\omega(h)\right)$$

$$= \left(1+\left[\Delta^2 h(0)\right]^{-1}\lambda_n 2^{-1}[\Delta^3 h(1)+\Delta h(2)]\right)\left(1-\lambda_n 6^{-1}\omega(h)\right)$$

$$= 1-\lambda 6^{-1}\left(\omega(h)-\left[\Delta^2 h(0)\right]^{-1}[\Delta^3 h(1)+\Delta h(2)]+O(\lambda_n^2)\right)$$

$$= 1-\frac{\lambda_n}{6(\Delta^2 h(0))^2}\Big(\Delta^3 h^2(0)-6\Delta h(1)\Delta^2 h(0)-2\Delta^3 h(0)[h(1)+2\Delta h(0)]$$

$$-3\Delta^2 h(0)[\Delta^3 h(1)+\Delta h(2)]+O(\lambda_n^2)\Big) = 1-\frac{\lambda_n(\Delta^3 h(0))^2}{6(\Delta^2 h(0))^2}+O(\lambda_n^2). \qquad \square$$

**Proof of Lemma 2.3**. The case $d=0$ is presented by Ivchenko and Medvedev (1978). Let $d\neq 0$. Rewrite the PDS in the form (2.9), where

$$\bar{\psi}_d(x) = \frac{2}{d(d+1)}\left[(x/\lambda_n)^{d+1}-(d+1)(x/\lambda_n)+d\right]$$

$$= \frac{2}{d(d+1)}\left[\left(1+\frac{x-\lambda_n}{\lambda_n}\right)^{1+d}-(1+d)\left(1+\frac{x-\lambda_n}{\lambda_n}\right)+d\right].$$

We have $\bar{\psi}_d(\lambda_n)=0$, $\bar{\psi}'_d(\lambda_n)=0$, $\bar{\psi}''_d(\lambda_n)=2\lambda_n^{-2}$, $\bar{\psi}'''_d(\lambda_n)=2(d-1)\lambda_n^{-3}$, $\bar{\psi}_d^{(4)}(\lambda_n)=2(d-1)(d-2)\lambda_n^{-4}$. Set $\bar{\xi}_n = (\xi-\lambda_n)/\lambda_n$, then $\bar{\xi}_n = O_p(\lambda_n^{-1/2})$, since $\sqrt{\lambda_n}\bar{\xi}_n$ has an asymptotically normal distribution as $\lambda_n \to \infty$. Apply these facts in Taylor expansion formula for the $\bar{\psi}_d(x)$ to get

$$\bar{\psi}_d(\xi) = \bar{\xi}_n^2 + \frac{d-1}{3}\bar{\xi}_n^3 + \frac{(d-1)(d-2)}{12}\bar{\xi}_n^4 + O_p(\bar{\xi}_n^5). \tag{3.34}$$

Put $\mu_l = E(\xi-\lambda_n)^l$. Then $\mu_2 = \mu_3 = \lambda_n$, $\mu_4 = 3\lambda_n^2+\lambda_n$, $\mu_5 = 10\lambda_n^2+\lambda_n$, $\mu_6 = 15\lambda_n^3+25\lambda_n^2+\lambda_n$, $\mu_7 = O(\lambda_n^3)$. Using these facts and (3.34) we obtain

$$E\bar{\psi}_d(\xi) = \frac{1}{\lambda_n} + \frac{d-1}{3}\frac{1}{\lambda_n^2} + \frac{(d-1)(d-2)}{12}(\frac{3}{\lambda_n^2} + \frac{1}{\lambda_n^3}) + O(\lambda_n^{-3}), \quad (3.35)$$

$$(E\bar{\psi}_d(\xi))^2 = \frac{1}{\lambda_n^2} + \frac{1}{\lambda_n^3}(\frac{2(d-1)}{3} + \frac{(d-1)(d-2)}{2}) + O(\lambda_n^{-4}),$$

$$E\bar{\psi}_d^2(\xi) = E\bar{\xi}_n^4 + 2\frac{d-1}{3}E\bar{\xi}_n^5 + (\frac{(d-1)^2}{9} + \frac{2(d-1)(d-2)}{12})E\bar{\xi}_n^6 + O(\lambda_n^{-4})$$

$$= \frac{3}{\lambda_n^2} + \frac{1}{\lambda_n^3}(1 + \frac{20(d-1)}{3} + \frac{(d-1)^2}{9} + \frac{5(d-1)(d-2)}{2}),$$

$$\tau_n = \frac{1}{\lambda_n}(E\bar{\xi}_n^2(\xi - \lambda_n) + \frac{d-1}{3}E\bar{\xi}_n^3(\xi - \lambda_n) + EO_p(\bar{\xi}_n^4(\xi - \lambda_n)))$$

$$= \frac{d}{\lambda_n^2} + \frac{1}{\lambda_n^3}(\frac{d-1}{3} + \frac{5(d-1)(d-2)}{6}) + O(\lambda_n^{-4}) = \frac{d}{\lambda_n^2}(1 + O(\lambda_n^{-1})),$$

$$\lambda_n \tau_n^2 = \frac{d^2}{\lambda_n^3}(1 + O(\lambda_n^{-1})).$$

Thus

$$\sigma_N^2(\bar{\psi}_d) = E\bar{\psi}_d^2(\xi) - (E\bar{\psi}_d(\xi))^2 - \lambda_n \tau_n^2 = \frac{2}{\lambda_n^2}\left(1 + \frac{1}{2\lambda_n}(d^2 - 1 + \frac{(d-1)^2}{9}) + O(\lambda_n^{-2})\right). \quad (3.36)$$

Further, by (3.34)

$$E(\bar{\psi}_d(\xi)(\xi - \lambda_n)^2 = E\bar{\xi}_n^2(\xi - \lambda_n)^2 + \frac{d-1}{3}E\bar{\xi}_n^3(\xi - \lambda_n)^2$$

$$+ \frac{(d-1)(d-2)}{12}E\bar{\xi}_n^4(\xi - \lambda_n)^2 + O(\lambda_n^{-2})$$

$$= 3 + \frac{1}{\lambda_n}\left(1 + \frac{10(d-1)}{3} + \frac{15(d-1)(d-2)}{12}\right) + O(\frac{1}{\lambda_n^2}),$$

$$E(\bar{\psi}_d(\xi)(\xi - \lambda_n)) = \frac{1}{\lambda_n} + \frac{d-1}{3}(\frac{3}{\lambda_n} + \frac{1}{\lambda_n^2}) + O(\lambda_n^{-2}) = \frac{d}{\lambda_n}\left(1 + O\left(\frac{1}{\lambda_n}\right)\right).$$

Applying these equalities and (3.35) we obtain

$$\text{cov}\left(h(\xi) - \tau_n \xi, \xi^2 - (2\lambda_n + 1)\xi\right) = E(\bar{\psi}_d(\xi)(\xi - \lambda_n)^2 - E\bar{\psi}_d(\xi)(\xi - \lambda_n) - \lambda_n E\bar{\psi}_d(\xi)$$

$$= 2 + \frac{1}{\lambda_n}\left(1 + \frac{10(d-1)}{3} + \frac{15(d-1)(d-2)}{12} - d - \frac{d-1}{3} - \frac{(d-1)(d-2)}{4}\right) + O(\lambda_n^{-2})$$

$$= 2 + \frac{1}{\lambda_n}(d^2 - d) + O(\lambda_n^{-2}). \quad (3.37)$$

Now by (3.2), (3.36) and (3.37) we obtain

$$\rho(h, \lambda_n) = \frac{2 + \frac{1}{\lambda_n}(d^2 - d) + O(\lambda_n^{-2})}{2\sqrt{1 + (2\lambda_n)^{-1}(d^2 - 1 + 9^{-1}(d-1)^2) + O(\lambda_n^{-2})}}$$

$$= \left[1 + (2\lambda_n)^{-1}(d^2 - d) + O(\lambda_n^{-2})\right]\left[1 - (2\lambda_n)^{-1}(d^2 - 1 + 9^{-1}(d-1)^2) + O(\lambda_n^{-2})\right]$$

$$= 1 - \frac{(d-1)^2}{6\lambda_n} + O(\lambda_n^{-2}). \qquad \square$$

**Conclusion**

The concept of an intermediate (between the Pitman and Bahadur settings) asymptotic relative efficiency (IARE) of tests of several types in the classical problem of testing uniformity [0,1] has been elaborated in the literature. We applied this setting to study of IARE of symmetric tests in discrete model, namely in the problem of testing of uniformity of a multinomial distribution. A formula for the IARE of arbitrary two symmetric tests, satisfying some conditions, is derived. A key point here is the validity of the large deviation results for the statistics of the comparing tests. The formula depends on the asymptotic correlation coefficients, counted under the hypothesis, of each test statistics with the chi-square statistic, and on the asymptotic behavior of $N$ -the number of cells, which assumed to be taken as a positive function of the continuous variable $x$ and has a power tail of index $q \in (0,2)$, i.e., $N(x) \sim cx^q$ as $x \to \infty$. The general result is applied to study IARE of the chi-squared test wrt symmetric tests for: sparse models, where $q=1$, very sparse models, where it is assumed $1 < q \leq 4/3$, and dense models, where $0 < q < 1$. The proven results and still unsolved problems are as follows.

*Sparse models*. It is achieved conflicting results: The chi-square test is superior within class of $h$-tests for the alternatives which converge to the hypothesis at the rate $\varepsilon(n) = o(n^{-1/3})$, as it is in the Pitman sense, see Holst (1972), Ivchenko and Medvedev (1978). But the IARE of the chi-square test wrt, for instance, tests based on the log-likelihood statistic, the Freeman-Tukey statistic, CS (1.5) can't be positive for the intermediate alternatives which far from the hypothesis at the distance $\varepsilon(n) \gg n^{-1/3} \log^{2/3} n$, as it is in the Bahadur efficiency of the chi-square test wrt the log-likelihood test, see Quine and Robinson (1985).

IARE of chi-square test wrt $h$-tests for the alternatives such that $c_1 n^{-1/3} \leq \varepsilon(n) \leq c_2 n^{-1/3} \log^{2/3} n$ remains unsolved.

*Very sparse models*. IARE of chi-square test wrt $\psi_d$-tests for such models depend on the decreasing rate of $\lambda_n$, the tuning parameter $d$ and considering family of alternatives. In particularly, if $d \in (-1, 1/2)$ and $\lambda_n \gg n^{-(1-2d^*)/6}$ then again conflict results are achieved: $\psi_d$-test and chi-square test are equally efficient for alternatives such that $\varepsilon(n) \leq n^{-1/3}$, but for the alternatives specified by condition $(n\lambda_n)^{-1/3} \log^{2/3}(N^2/n) \ll \varepsilon(n) \ll (n\lambda_n^2)^{-1/4}$ IARE of chi-square test wrt $\psi_d$-test is equal to zero. For the case $d \geq 1/2$ it is shown that chi-square test and $\psi_d$-test and are equally efficient for the alternatives specified by $\varepsilon(n) \ll (\lambda_n^{1-d}/n^d)^{1/(2d+1)}$. We emphasize that $n\lambda_n^3 \to \infty$ (i.e., $n \ll N \ll n^{4/3}$) is a necessary condition for a family of alternatives with $\varepsilon(n) \ll (\lambda_n^{1-d}/n^d)^{1/(2d+1)}$ to be a non-empty sub-family of $\Im_{alt}$. Note that $(\lambda_n^{1-d}/n^d)^{1/(2d+1)} \ll n^{-1/3}$, if $d \geq 1/2$.

IARE of chi-square test wrt $\psi_d$-tests for the alternatives such that $n^{-1/3} \leq \varepsilon(n) \leq (n\lambda_n)^{-1/3}$ or $\varepsilon(n) \geq (n\lambda_n^2)^{-1/4}$ if $d \in (-1, 1/2)$, and for alternatives such that $\varepsilon(n) \geq (\lambda_n^{1-d}/n^d)^{1/(2d+1)}$ if $d > 1$ remain open.

*Dense models.* It is pointed out that in this case under the hypothesis each PDS asymptotically in distribution coincides with the chi-square statistics. This fact seems is reflected to IARE of chi-square test wrt $\psi_d$-test, specifically the chi-square test and $\psi_d$ test are equally efficient for alternatives such that $\varepsilon(n) \ll (n\lambda_n^2)^{-1/3}$ and each $d > -1$ (here no condition to increasing rate of $N$ is assumed); the same verdict is valid for every alternative of $\Im_{alt}$ and $d = 0$ if $N = o(n^{3/8})$. These extend to the intermediate setting the known result that: if $\lambda_n \to \infty$ then the Pitman efficiency of the chi-square test wrt the $\psi_d$-test is equal to 1. However, I conjecture that aforesaid property of PDS effects to IARE properties of $h$-tests if $N \ll \sqrt{n}$. This is confirmed by the fact that the log-likelihood and chi-square tests have the same asymptotic efficiency in term of $\alpha$ - IAE for alternatives such that $\varepsilon(n) \ll (n\lambda_n^2)^{-1/4}$ if $N \ll \sqrt{n}$ (note that alternatives with $\varepsilon(n) \geq (n\lambda_n^2)^{-1/4}$ are out of the family $\Im_{alt}$ if $N \ll \sqrt{n}$), see Mirakhmedov (2021). But (ibid) $\alpha$ - IARE of chi-square test is much inferior wrt tests satisfying the Cramer condition, for instance the log-likelihood and Freeman-Tukey tests, if $N \geq \sqrt{n}$ and alternatives are such that $\varepsilon(n) \gg (n\lambda_n)^{-1/3} \log^{2/3}(N^2/n)$; this family may contain the alternatives with $\varepsilon(n) \geq (n\lambda_n^2)^{-1/4}$.

Thus, IARE of the chi-square test wrt $h$-tests in the case $\lambda_n \to \infty$ but $\lambda_n \ll N$ completely, and in the case $\lambda_n \gg N$ for alternatives such that $(n\lambda_n)^{-1/3} \leq \varepsilon(n) \ll (n\lambda_n^2)^{-1/4}$ remains open.

Progress in the study of IARE $h$-tests depends on the progress of the results on the probability of large deviations of symmetric statistics (1.2).

**Acknowledgements.** The author thanks the anonymous reviewer for his very friendly and valuable constructive comments that lead to a substantial improved version of the paper.

---

[i] Mirakhmedov S .M. is the former Mirakhmedov S .A.